\documentclass{ws-procs9x6}

\newcommand\R{{\mathbb R}}

\newcommand\C{{\mathbb C}}

\def\AA{{\mathcal A}}
\def\BB{{\mathcal B}}

\def\HH{{\mathcal H}}

\def\LL{{\mathcal L}}

\def\RR{{\mathcal R}}

\def\TT{{\mathcal T}}

\newtheorem{theo}{Theorem}

\newtheorem{cor}[theo]{Corollary}
\newtheorem{rem}[theo]{Remark}

\renewcommand{\theequation}{\thesection.\arabic{equation}}

\newcommand{\beqn}{\begin{equation}}
\newcommand{\eeqn}{\end{equation}}
\newcommand{\bear}{\begin{eqnarray}}
\newcommand{\eear}{\end{eqnarray}}
\newcommand{\bean}{\begin{eqnarray*}}
\newcommand{\eean}{\end{eqnarray*}}

\begin{document}

\title{Enlarging the functional space of decay estimates on
  semigroups}

\author{C. MOUHOT$^*$}

\address{CNRS \& DMA, \'ENS Paris,\\
45, rue d'Ulm, F-75230 Paris cedex 05\\
$^*$E-mail: Clement.Mouhot@ens.fr}



\begin{abstract}
  This note briefly presents a new method for enlarging the functional
  space of a ``spectral-gap-like'' estimate of exponential decay on a
  semigroup. A particular case of the method was first devised in
  Ref.~\refcite{Mcmp} for the spatially homogeneous Boltzmann equation, and
  a variant was used in Ref.~\refcite{MMcmp} in the same context for
  inelastic collisions. We present a generalized abstract version of it, a
  short proof of the algebraic core of the method, and a new application to
  the Fokker-Planck equation. More details and other applications shall be
  found in the work in preparation Ref.~\refcite{GMM} (another application
  to quantum kinetic theory can be found in the work in
  preparation Ref.~\refcite{many}).
\end{abstract}


\bodymatter

\section{The ``space enlargement'' issue}
\label{sec:space-enlarg-issue}

Consider a Hilbert space $\HH$, a (possibly unbounded) linear operator
$\TT$ on $\HH$ which generates a strongly continuous semigroup $e^{t
  \, \TT}$ with spectrum $\Sigma(\TT)$. Assume that for some Hilbert
{\em subspace} $H \subset \HH$ the restricted operator $T := \TT
\big|_H$ generates a strongly continuous semigroup $e^{t \, T}$ with
spectrum $\Sigma(T)$ in $H$.

Assume some ``spectral-gap-like'' information on $\Sigma(T)$,
typically when $T$ is self-adjoint assume
$$
\forall \, f \in H, \ f \, \bot \, \mbox{Null}(T), \quad \left\| e^{t \, T}
  \, f \right\|_H \le e^{\lambda \,t} \, \| f \|_H, \quad \lambda <0.
$$

An important class of applications is the following: $\TT$ is a partial
differential operator (acting on a large class of function on $\R^d$, say
$L^1$), with equilibrium $\mu$ and detailed spectral information available
in a much smaller space $H = L^2(\mu^{-1})$ where it is symmetric. The
latter space is much smaller than $\HH$ in the sense that it requires a
stronger decay condition, e.g. when $\mu$ is a gaussian in statistical
mechanics.

The question addressed here is: can one deduce from the spectral-gap
information in the space $H$ some spectral-gap information in the larger
space $\HH$, and if possible in a quantitative way? More explicitly, does
$e^{t \, \TT}$ have the same decay property as $e^{t \, T}$ above?

We give a positive answer for a class of operators $\TT$ which split into a
part $\AA$ ``regularizing'' $\HH$ into $H$ and a coercive part $\BB$. We
then show that, under some assumption on the potential force, the
Fokker-Planck equation belongs to this class and, as a consequence, we
prove that its spectral gap property can be extended from the linearization
space (with gaussian decay) to larger $L^2$ spaces with, say, polynomial
weights.

\section{The abstract result}
\label{sec:abstract}
\setcounter{equation}{0}
\setcounter{theo}{0}

Let us start with an almost equivalent condition for the decay of the
semigroup in terms of a uniform bound on a vertical line for the
resolvent. We omit the proof to keep this note short. It can be found
in Ref.~\refcite{GMM} and it mainly relies on a careful use of the Parseval
identity between the resolvent operator and the semigroup.

For some closed densely defined unbounded operator $T$ in a Hilbert space
$E$, denote by by $R(z) = (T-z)^{-1}$, $z \not \in \Sigma(T)$ 
its resolvent operator, and $\LL(E)$ the space of bounded linear operators
on $E$. Finally for any $a \in \R$, define the half complex plane $\Delta_a
:= \{ z \in \C, \, \Re e \, z > a \}$.

\begin{theo}\label{AbstractTheo1} Assume for the operator $T$ in the
  Hilbert space $E$:

  \begin{itemize}
  \item[{\bf (H1)}] {\bf Localization of the spectrum}: $\Sigma (T)
    \subset (\overline{\Delta_a})^c \cup \{ \xi_1, \, ... \, , \xi_k
    \}$ with $a \in \R$, and $\xi_j \in \Delta_a$, $1\le j \le k$ some
    discrete eigenvalues;  
  \item[{\bf (H2)}] {\bf Control on the resolvent operators}: 
    $$
    \exists \, K >0, \quad \forall \, y \in \R, \,\,\, \quad \| R(a+i
    \, y) \|_{\LL(E)} \le K.
    $$
  \item[{\bf (H3)}] {\bf Weak control on the semigroup}: There exist $b,
    C_b \ge 0$ such that
    $$
    \forall \, t \ge 0 \qquad \| e^{t \, T} \|_{\LL(E)} \le C_b \, e^{b \,
      t}.
    $$ 
  \end{itemize}
  
  Then, for any $\lambda > a$, there exists $C_\lambda$ explicit from
  $a$, $b$, $C_b$, $K$ such that
  \begin{equation}\label{bddSlambda} 
    \forall \, t \ge 0, \quad 
    \left\| e^{t \, T} - \sum_{i=1}^k e^{\xi_i \, t} \, \Pi_i
    \right\|_{\LL(E)} 
    \le C_\lambda \, e^{\lambda \, t}  
  \end{equation}
  for the spectral projectors $\Pi_i$ of eigenvalues $\xi_i$.

  We also have the following converse result: assume
  \begin{equation}\label{bddSa}
    \forall \, t \ge 0, \quad 
    \left\| e^{t \, T} - \sum_{j=1}^k e^{\xi_j \, t} \, \Pi_j
    \right\|_{\LL(E)} \le C_a \, e^{a \, t}
  \end{equation}
  for some constants $a \in \R$, $C_a \in (0,\infty)$, some complex
  numbers $\xi_j \in \Delta_a$ and some operators $\Pi_j$ which
  all commute with $e^{t \, T}$. Then $T$ satisfies {\bf (H1)},
  {\bf (H2)}, {\bf (H3)}.
\end{theo}

\begin{rem}
  Assumption {\bf (H3)} is required in this theorem in order to obtain
  quantitative constants in the rate of decay. Therefore, under assumptions
  {\bf (H1)} and {\bf (H3)}, assertions {\bf (H2)} and
  Eq.~\eqref{bddSlambda} are equivalent in a quantitative way.
\end{rem}


The following theorem is the core of the method:

\begin{theo}\label{AbstractTheo3} 
  Assume that $\TT$ is a closed unbounded densely defined operator in a
  Hilbert space $\HH$, and that $T := \TT \big|_H$ is a closed unbounded
  densely defined operator in a Hilbert subspace $H \subset \HH$ which
  satisfies {\bf (H1)} and {\bf (H2)} (with $E=H$). Assume moreover that
  $\TT$ satisfies:
  \begin{itemize}
  \item[{\bf (H4)}] {\bf Decomposition:} $\TT = \AA + \BB$ where $\AA$ and
    $\BB$ are closed unbounded densely defined operators with domains
    included in the one of $\TT$ such that
    \begin{itemize}
    \item for some $r >0$, the operator $\BB-\xi$ is invertible with
      uniform bound for any $\xi \in \Delta_a \setminus ( \cup_{i=1} ^k
      B(\xi_i,r))$ (where every balls $B(\xi_i,r)$ are strictly included
      in $\Delta_a$);
    \item $B = \BB_{|H}$ is well-defined as a closed unbounded densely
      defined operator with domain included in the one of $T$, and $B-\xi$
      is invertible for any $\xi \in \Delta_a \setminus ( \cup_{i=1} ^k
      B(\xi_i,r))$;
    \item $\AA \, (\BB-\xi)^{-1} : \HH \to H$ and $(\BB-\xi)^{-1} \, \AA :
      \HH \to H$ are bounded for any $\xi \in \Delta_a \setminus (
      \cup_{i=1} ^k B(\xi_i,r))$.
    \end{itemize}
  \end{itemize}
  Then $\TT$ satisfies {\bf (H2)} in the space $E=\HH$ (with
  constructive bounds in terms of the above assumptions). 
\end{theo}

The proof of the following corollary is immediate by combining Theorem
\ref{AbstractTheo3} and Theorems \ref{AbstractTheo1}.

\begin{cor}
  Assume that $T$ satisfies {\bf (H1)} and {\bf (H2)} in the space $H$
  and $\TT$ satisfies {\bf (H3)} and {\bf (H4)} in the space $\HH$.
  Assume moreover that the eigenvalues of $\TT$ in $\Delta_a$ are the
  same as those of $T$, that is $\{\xi_1,\dots,\xi_k\}$. Then the
  conclusion of Theorem~\ref{AbstractTheo1} holds in the space $\HH$:
  for any $\lambda > a$, there exists $C_\lambda$ explicit from $a$,
  $b$, $C_b$, $K$ such that
  \begin{equation*}
    \forall \, t \ge 0, \quad 
    \left\| e^{t \, \TT} - \sum_{i=1}^k e^{\xi_i \, t} \, \Pi_i
    \right\|_{\LL(\HH)} 
    \le C_\lambda \, e^{\lambda \, t}.  
  \end{equation*}
\end{cor}

\begin{rem}
  If $r$ can be taken as small as wanted {\bf (H4)} (for some
  decompositions depending on $r$), it can be proved that the eigenvalues
  of $\TT$ in $\Delta_a$ are the same as those of $T$ in $\Delta_a$ (that
  is $\{\xi_1,\dots,\xi_k\}$) and this assumption can be relaxed.
\end{rem}

\begin{rem}
  Thanks to the reciprocal part of Theorem~\ref{AbstractTheo1},
  assumption {\bf (H2)} on $T$ can be replaced by assuming a decay on
  the semigroup:
  \begin{equation*} 
    \forall \, t \ge 0, \quad 
    \left\| e^{t \, T} - \sum_{j=1}^k e^{\xi_j \, t} \, \Pi_j
    \right\|_{L(H)} 
    \le C_\lambda \, e^{\lambda \, t}.  
  \end{equation*}
\end{rem}

\medskip\noindent {\sl Proof of Theorem~\ref{AbstractTheo3}.} 
Assume that $k=1$ and $\xi_1 = 0$ for the sake of simplicity, the
proof being similar in the general case.

Take $\xi \notin \Delta_a \backslash B(0,r)$ and define
$$
U(\xi) := \BB(\xi)^{-1} -  R(\xi) \, \AA \,\BB(\xi)^{-1},
$$
where $R(\xi)$ is the resolvent of $T$ in $H$ and $\BB(\xi) = \BB -
\xi$. Since by assumption $\BB(\xi)^{-1} : H \to H$, $\AA \,
\BB(\xi)^{-1} : \HH \to H$ and $R(\xi) : H \to H$ are bounded
operators, $U(\xi) : H \to H$ is well-defined and bounded from $\HH$ to
$\HH$. Then,
\begin{eqnarray*} 
  (\TT-\xi) \, U(\xi)
  &=& (\AA+\BB(\xi)) \, \BB(\xi)^{-1} 
  - (\TT-\xi) \, R(\xi) \, \AA \, \BB(\xi)^{-1} \\
  &=& \AA \, \BB(\xi)^{-1} + \mbox{Id}_\HH 
  - (\TT-\xi)  \,  R(\xi)  \, \AA \, \BB(\xi)^{-1} \\
  &=& \AA \, \BB(\xi)^{-1} + \mbox{Id}_\HH - \AA \, \BB(\xi)^{-1} =
  \mbox{Id}_\HH.
\end{eqnarray*} 
To be more precise, introduce the canonical injection $J : H \to \HH$ and
use that $R = J \, R$, $\AA = J \, \AA$, $\TT \, J = J \, T$ to write:
\begin{eqnarray*}
  (\TT-\xi) \,
  R(\xi) \, \AA \, \BB(\xi)^{-1}
  &=& (\TT-\xi) \, J \, R(\xi) \, \AA \, \BB(\xi)^{-1} = 
  J \, (T-\xi) \, R(\xi) \, \AA \, \BB(\xi)^{-1} \\
  &=& J \, \mbox{Id}_H \, \AA \, \BB(\xi)^{-1} = J \, \AA \,
  \BB(\xi)^{-1} = \AA \, \BB(\xi)^{-1}.
\end{eqnarray*} 
The operator $\TT - \xi$ is also one-to-one. Indeed, if $g$ satisfies
$$
g \in \mbox{Dom}(\TT), \qquad   (\TT - \xi) \, g = 0,
$$
the decomposition {\bf (H4)} yields
$$
\BB (\xi) \, g = - \AA \, g \in H,
$$
and therefore $g \in \mbox{Domain}(B) \subset \mbox{Domain}(T) \subset H$
because $B(\xi) = \BB(\xi)_{|H}$ is invertible on $H$. We conclude
that $g = 0$ since $T-\xi$ is one-to-one.  As a conclusion, $U(\xi)$
is the inverse of $\TT - \xi$ which in turn implies that $\xi \notin
\Sigma(\TT)$ and $\RR(\xi) = U(\xi)$ satisfies the announced
estimate. This concludes the proof.

\section{Application to the Fokker-Planck equation}
\label{sec:appl-fokk-planck}

In this section we are concerned with the Fokker-Planck equation
$$
\partial_t f = \TT \, f := \hbox{div} (\nabla f + E \, f)
$$
for the real valued density function $f = f(t,x)$, $t \ge 0$, $x \in
\R^d$. In this equation $E = E(x) \in \R^d$ is a given force field, written
as
\begin{equation}\label{FP-H1-1}
E = \nabla U + F
\end{equation}
where the potential $U : \R^d \to \R$ is such that $\mu (dx) =
e^{-U(x)} \, dx$ is a probability measure satisfying the ``Poincar\'e
inequality condition'': there exists $\lambda_P< 0$ such that
\begin{equation}\label{galineqPoincare}
  - \int |\nabla u |^2 \, \mu(dv) \le 
  \lambda_P \int u^2 \, \mu(dv) \qquad 
  \forall \, u \in H^1(\R^d), \,\,\, \int u \,  d\mu = 0,
\end{equation}
({\it cf.} for instance Refs. \refcite{deuschstroock,ledouxlivre,be} and
the references therein).  The additionnal force field $F$ satisfies
\begin{equation}\label{FP-H1-2}
  \nabla \cdot F = 0, \qquad  \nabla U \cdot F = 0, 
  \qquad |F| \le C(1 + |\nabla U|).
\end{equation}


Thanks to that structural assumptions we can split $\TT$ between a
symmetric term and a skew-symmetric term:
$$
\TT = \TT^{s} + \TT^{as}, \quad \TT^s f = \hbox{div} (\nabla f + \nabla
U \, f), \quad \TT^{as} f= \hbox{div} \, (F \, f).
$$
The operator $\TT^s$ is symmetric in $H = L^2(\mu^{-1})$
$$
\langle \TT^s f, g \rangle_H = - \int \nabla(f/\mu) \cdot
\nabla(g/\mu) \, \mu = \langle  f, \TT^s g \rangle_H,
$$
while the operator $\TT^{as}$ is anti-symmetric in $\HH = L^2(m^{-1})$
for any weight function $m^{-1} (v) = \theta (U(x))$, with $\theta :
\R_+ \to \R_+$:
\begin{eqnarray*}
(\TT^{as} f, g)_\HH &=& \int \Big[  
(\nabla \cdot F) \, f + F \cdot \nabla f \Big]  \, g \, m^{-1} 
= - \int  f  \, \nabla \cdot ( F \, g \, m^{-1}) \\
&=& - \int  f  \, \Big[ 
(\nabla \cdot F ) \, g \, m^{-1} +( F \cdot \nabla g) \, m^{-1} + 
g \, \theta'(U) \, (\nabla U \cdot F)\Big] \\
&=& - \int \Big[ (\nabla \cdot F ) \, g + F \cdot \nabla g \Big] \, f \,
m^{-1} = - \langle f, \TT^{as} g \rangle_\HH. 
\end{eqnarray*}
As an important consequence, we have
$$
\langle \TT^{as} f, f \rangle _\HH = \int \nabla \cdot (F \, f) \, f
\, \theta(U) = 0.
$$

In $H := L^2(\mu^{-1})$ the restricted operator $T := \TT \big|_H$ is
non-positive, its first eigenvalue is $0$ associated to the eigenspace $\R
\mu$, and it has a spectral gap thanks to the Poincar\'e inequality:
$$
\int f \, (Tf) \, \mu^{-1} =- \int \mu \, |\nabla(f/\mu))|^2 \le
\lambda_2 \, \| f - \langle f \rangle \|_H^2.
$$


A natural question to ask is whether it is possible to obtain an
exponential decay on the semigroup in a space larger than $H$. The
following theorem gives an answer in $L^2$ spaces with polynomial or
``stretched'' exponential weights. The proof follows from the application
of the abstract method and some careful computations on the Dirichlet form
in the larger space.

 \begin{theo}\label{theomFP2} 
   Let $\mu = e^{-U}$ with $U(v) = (1+|x|^2)^{s/2}$, $s \ge 1$ (so that
   Poincar\'e inequality holds for $\mu$).
   
   Let $m \in C^2(\R^d)$ be a weight function such that $m^{-1}(x) =
   \theta(U(x))$ with 
   $\theta(x) = (1+|x|^2)^{k/2}$ with $k > d$ or $\theta (x) =
   e^{(1+|x|^2)^{k/2}}$ with $k \in (0,1)$. Let us define $\HH :=
   L^2(m^{-1})$.

   Then there exist explicit $\lambda \in (-\lambda_P,0)$ and $C_\lambda
   \in [1,\infty)$ such that
   $$
   \forall \, f_0 \in \HH, \qquad \forall \, t \ge 0 \quad \| f_t
   - \langle f_0 \rangle \, \mu \|_\HH \le C_\lambda \, e^{\lambda
     \, t} \, \| f_0 - \langle f_0 \rangle \, \mu \|_\HH.
   $$
\end{theo}

\begin{rem}\label{remFP1} 
  In this theorem $C_\lambda > 1$ is allowed, which means that we do
  not prove that the Dirichlet form of $\TT$ has a sign.
  \end{rem}
  
\begin{rem}
  The smoothness assumption on $U$ and $m$ at the origin can be
  relaxed.
\end{rem}

\noindent {\bf Acknowledgements.} The author thanks St\'ephane Mischler for
useful discussion during the preparation of this note. He also wishes to
thank the Award No. KUK-I1-007-43, funded by the King Abdullah University
of Science and Technology (KAUST) for the funding provided for his repeated
visits at Cambridge University.


\bibliographystyle{ws-procs9x6}
\bibliography{Proc-Isaac-Mouhot}

\end{document}